\documentclass[12pt]{amsart}
\usepackage{amssymb,amsmath,amscd}
\newtheorem{thm}{Theorem}[section]
\newtheorem{corollary}[thm]{Corollary}
\newtheorem{prop}[thm]{Proposition}
\newtheorem{lemma}[thm]{Lemma}
\newtheorem{fact}[thm]{Fact}

\theoremstyle{definition}
\newtheorem{defn}[thm]{Definition}
\newtheorem{example}[thm]{Example}
\newtheorem{remark}[thm]{Remark}

\newcommand{\bt}{\begin{thm}}
\newcommand{\et}{\end{thm}}
\newcommand{\bp}{\begin{prop}}
\newcommand{\ep}{\end{prop}}
\newcommand{\bd}{\begin{defn}}
\newcommand{\ed}{\end{defn}}
\newcommand{\bl}{\begin{lemma}}
\newcommand{\el}{\end{lemma}}
\newcommand{\bfa}{\begin{fact}}
\newcommand{\efa}{\end{fact}}
\newcommand{\bc}{\begin{corollary}}
\newcommand{\ec}{\end{corollary}}
\newcommand{\bex}{\begin{example}}
\newcommand{\eex}{\end{example}}
\newcommand{\br}{\begin{remark}}
\newcommand{\er}{\end{remark}}

\newcommand{\isom}{\cong}

\newcommand{\Hom}{{\mathcal H}om}

\newcommand{\Pic}{{\mbox{Pic}}}
\newcommand{\Supp}{{\mbox{Supp}}}
\newcommand{\Hilb}{{\mbox{Hilb}}}

\newcommand{\Ext}{{\mathcal E}xt}

\newcommand{\sotto}[2]{#1_{#2}}

\newcommand{\rrr}{\rightarrow}
\newcommand{\ra}{\rightarrow}

\newcommand{\ideal}[1]{\sotto {{\mathcal I}}{#1}}

\newcommand{\exact}[3]
{0 \rrr #1 \rrr #2
\rrr #3 \rrr 0}

\newcommand{\pso}{\mathbb{P}^3}

\newcommand{\PP}{\mathbb{P}}
\newcommand{\pp}{\mathbb{P}}

\newcommand{\Ptwo}{\mathbb{P}^2}

\newcommand{\Cpl}{\mathbb{C}}
\newcommand{\Z}{\mathbb{Z}}

\newcommand{\cnn}{{\mathcal N}}
\newcommand{\coo}{{\mathcal O}}

\newcommand{\caf}{{\mathcal F}}

\newcommand{\mathcall}{{\mathcal L}}
\newcommand{\scrl}{{\mathcal L}}
\newcommand{\cae}{{\mathcal E}}
\newcommand{\cag}{{\mathcal G}}
\newcommand{\cah}{{\mathcal H}}

\begin{document}

\title{Curves on a Double Surface}

\author{Scott Nollet}
\address{Department of Mathematics,
Texas Christian University, Fort Worth, TX 76129, USA}
\email{s.nollet@tcu.edu}

\author{Enrico Schlesinger}
\address{Dipartimento di Matematica, Politecnico di Milano, Piazza
Leonardo da Vinci 32,
20133 Milano, Italy}
\email{enrsch@mate.polimi.it}

\date{}
\dedicatory{Dedicated to Silvio Greco on the occasion of his sixtieth birthday}

\begin{abstract}
     Let $X$ be a doubling of a smooth surface $F$ in a smooth
     threefold and let $C \subset X$ be a locally Cohen-Macaulay curve.
     Then $C$ gives rise to two effective divisors on $F$, namely the
     curve part $P$ of $C \cap F$ and the curve $R$ residual
     to $C \cap F$ in $C$.
     We show that a general deformation of $R$ on $F$ lifts to a
     deformation of $C$ on $X$ when a certain cohomology group
     vanishes and give applications to the study of
     Hilbert schemes of locally Cohen-Macaulay space curves.
\end{abstract}


\maketitle

\section{Introduction}

It is usually difficult to determine when a fixed curve $C \subset \pso$
is in the closure of another family of curves.
Beyond semicontinuity conditions, there are few known obstructions.
Hartshorne showed that curves of certain degree and genus cannot
specialize to stick figures by analyzing the specific quadric and
cubic surfaces on which the curves lie \cite{zeuthen}.
The problem reduces to that of linear equivalence if all the curves lie
on a {\it fixed smooth} surface $F$, since the families of divisor classes
are both open and closed in the corresponding Hilbert scheme \cite{F}.
When the curves in question are nonreduced, smooth surfaces are of little help.

In our analysis of curves of degree four in $\pso$ \cite{NS}, we
used families of curves on double planes and double quadric surfaces to
produce various specializations in the Hilbert schemes: these were
critical in determining irreducible components and showing
connectedness. If $X$ is a doubling of a smooth surface $F$
in a smooth threefold - or more generally, if $X$ is a ribbon
supported on $F$ in the sense of Bayer and Eisenbud \cite{beis} -
we will describe the Hilbert scheme of curves on $X$, using
as our model the rather complete study of curves on a double plane in $\pso$
\cite{2h}.

In section 2 we describe the natural triple $T(C)$ associated to a curve
$C \subset X$ defined in \cite{2h}: the scheme-theoretic
intersection $C \cap F$ has a divisorial part $P$ and a zero-dimensional
part $Z$, and when we form the residual curve $R$ to $C \cap F$ in $C$, we
obtain the triple $T(C)=\{Z,R,P\}$. Here $Z$ is a generalized Gorenstein
divisor on $R$ and $R \subset P$ are effective divisors on $F$.
We describe the curves $C$ giving rise to a fixed triple $\{Z,R,P\}$
and give practical conditions (\ref{exist} and \ref{practical})
under which this space is non-empty. The existence of such curves $C$
is subtle when our conditions fail.

Using the triples above, we stratify $H_{d,g}(X)$ in section 3 to obtain
locally closed $H_{z,r,p} \subset H_{d,g} (X)$ with natural projection
maps $t$ to the relevant Hilbert flag schemes $D_{z,r,p}$.
Relativizing results from section 2, we
find (\ref{structureofpi}) that $t$ has the local structure of an open
immersion followed by an affine bundle projection over the locus
$V \subset D_{z,r,p}$ of triples satisfying $H^{1}(\coo_{R}(Z+P-F))=0$
(the fibres are nonempty if condition (2) of (\ref{exist}) holds).
Combining with the structure of the projection map
$D_{z,r} \ra H_{r}$ (\ref{dominant}), we find (Theorem \ref{limits}) that
if $C \subset X$ is a curve whose triple $\{Z,R,P\}$
satisfies $H^{1}(\coo_{R}(Z+P-F))=0$, then a general deformation
of $R$ lifts to a deformation of $C$. We close with some applications
to families of space curves.

\section{Curves on a ribbon} \label{two}

Let $F$ be a smooth surface  over an algebraically closed ground field $k$.
If $F$ is contained in a smooth threefold $T$ and $X$ is the effective divisor
$2F$ on $T$, then

\begin{enumerate}
\item
$\Supp X = F$;
\item
$\ideal{F,X} \cong \coo_{F} (-F)$ is an invertible $\coo_{F}$-module.
\end{enumerate}
In other words, $X$ is a {\it ribbon} over $F$
in the sense of Bayer and Eisenbud \cite{beis}.
Since $F$ is smooth, any ribbon $X$ is locally split, hence appears
locally as a doubling of $F$ in a smooth threefold.
We use the notation $\coo_{F} (-F)= \ideal{F}$ and
$\coo_F (F) = \cah om_{\coo_F} ( \ideal{F}, \coo_F)$.
Here we will further assume $X$ is projective, although many of
our constructions work more generally.

We will study curves on a ribbon $X$ over $F$ using the triples introduced
in ~\cite{2h}. We adopt the following conventions:
A subscheme $C \subset X$ is a {\it curve} if
all of its associated points have dimension one, thus $C$ is locally
Cohen-Macaulay of pure dimension one or empty. If $Y$ is a subscheme of
$X$, $\ideal{Y}$ denotes the ideal sheaf of $Y$ in $X$.
If  $R$  is a  Gorenstein  scheme and  $Z$  a  generalized divisor  on
$R$ \cite{hgd},  then $\coo_R  (Z)  =  \cah  om  (
\ideal{Z,R},  \coo_R)$  denotes the  reflexive sheaf  associated  to  the
divisor  $Z$. If further $R \subset  F$, we write  $\coo_R (Z-F)$  for
$\coo_R (Z) \otimes \coo_F (-F)$.

\bp \label{1.1}
To each curve $C$ in $X$ is associated a triple  $T(C) = \{Z,R,P \}$
in which $R \subset P$ are effective divisors on $F$, $Z \subset R$ is
Gorenstein and zero-dimensional (possibly empty), and
$$
\ideal{P, C} \cong \coo_R (Z-F).
$$
The arithmetic genera are related by
\begin{equation}
p_{a}(C)= p_{a}(P) + p_{a}(R) + \deg_R \coo_{R}(F) - \deg (Z) -1.
\end{equation}
\ep

\begin{proof}
We proceed as in~\cite[$\S 2$]{2h}. Extracting the
possible embedded points from the one dimensional scheme-theoretic
intersection $C \cap F \subset F$, we may write
$$
\ideal{C \cap F,F} = \ideal{Z,F} (-P)
$$
where $P$ is an effective divisor and $Z$ is zero-dimensional.
The inclusion $P \subset C \cap F$ yields a commutative diagram
\begin{equation} \label{diagram}\begin{array}{lllllllll}
          & & 0 & & 0 & & 0 & & \\
          & & \downarrow & & \downarrow & & \downarrow & & \\
          0 & \ra & \ideal{R,F}(-F) & \ra & \ideal{C,X} & \ra &
          \ideal{Z,F}(-P) & \ra & 0 \\
          & & \downarrow & & \downarrow & & \downarrow & & \\
          0 & \ra & \coo_F (-F) & \ra & \ideal{P,X} & \ra &
          \coo_{F}(-P) & \ra & 0 \\
          & & \downarrow & & \downarrow 
          & & \downarrow & & \\
          0 & \ra & \coo_{R}(-F) & \stackrel{\sigma}{\ra} & \mathcall &
\ra & \coo_{Z}(-P) & \ra & 0 \\
          & & \downarrow & & \downarrow & & \downarrow & & \\
          & & 0 & & 0 & & 0 & &  \\
\end{array}\end{equation}
which defines the residual scheme  $R$ to  $C \cap F$ in $C$.
The  inclusion $\coo_R  (-F) \hookrightarrow \coo_C$  shows that the
associated points of $R$ are among those of $C$, hence $R$ is a curve.
By construction, $P$ is the largest curve in $F \cap C$,
hence $R \subseteq P$ and $C \subset F$ if and only if $R$ is empty.


We now show that $Z$ is Gorenstein on $R$ and that
$\mathcall \cong \coo_{R}(Z-F)$ is a rank one reflexive
$\coo_{R}$-module.
In view of the bottom row of diagram
\ref{diagram}, the submodule $\ideal{R} \scrl \subset \scrl$ is
supported on $Z$, but $\scrl = \ideal{P,C} \subset \coo_{C}$ has
only associated points of dimension one (because $C$ is purely
one-dimensional), hence $\ideal{R} \scrl =0$ and $\scrl$ is an
$\coo_{R}$-module.
It follows that $\coo_{Z}(-P)$ is an $\coo_{R}$-module as well,
hence $Z \subset R$.

Applying the bifunctor $\Hom_{\coo_R} (\: - \:, \: - \: )$ to
the sequence $\exact{\ideal{Z,R}}{\coo_R}{\coo_Z}$  and
the bottom row of diagram~\ref{diagram} we obtain
\begin{equation*} \begin{array}{cccccccc}
            & & 0 & & 0 & &  \Hom_{\coo_R} (\coo_Z, \coo_Z (-P)) \\
            & & \downarrow & & \downarrow && \downarrow 
\mbox{\scriptsize{$\alpha$}}     \\
            0 & \ra &
  \coo_R (-F)  & \stackrel{\pi}{\ra} &
  \coo_R (Z-F) & \ra &
\Ext^{1}_{\coo_R} (\coo_Z, \coo_R (-F)) \\
            & & \downarrow & & \downarrow \mbox{\scriptsize{$\beta$}} & &
       \downarrow 0\\
            0 &        \ra &
  \scrl & \stackrel{\phi}{\ra} &
\Hom_{\coo_R} (\ideal{Z,R}, \scrl)          & \stackrel{\gamma}{\ra} &
\Ext^{1}_{\coo_R} (\coo_Z, \scrl)  \\
           & & \downarrow & & \downarrow & & \\
\coo_Z (-P) & \stackrel{\delta}{\ra} &
\coo_Z (-P) & \stackrel{0}{\ra} &
\Hom_{\coo_R} (\ideal{Z,R}, \coo_Z (-P)) &  &  \\
\end{array}\end{equation*}
The morphisms $\pi$, $\phi$ and $\alpha$ are injective because
$\Hom_{\coo_R} (\coo_Z, \scrl)=0$ as $\scrl$ has no zero dimensional associated
point.

Since $F$ is smooth, $R$ is Gorenstein, hence
$\omega_{Z} \cong \Ext^{1}_{\coo_R} (\coo_Z,\coo_R (-F))$ and $\alpha$
can be thought as a morphism $\coo_Z \ra \omega_Z$. Since $\coo_Z$ and
$\omega_Z$ have the same length, $\alpha$ is an isomorphism (which
explains the $0$ map at the right of the diagram) and $Z$ is
Gorenstein.

Thinking of $\scrl$ and $\coo_R (Z-F)$ as subsheaves of
$\Hom_{\coo_R} (\ideal{Z,R}, \scrl)$, the $0$ at the bottom
of the diagram yields $\scrl \subset \coo_R (Z-F)$
while the $0$ at the right gives $\coo_R (Z-F) \subset \scrl$, hence
$\scrl=\coo_R (Z-F)$.

For the arithmetic genus formula, note that $p_{a}(C)-p_{a}(P) =
- \chi \ideal{P,C} = - \chi \scrl$, which can be read off from the
bottom row of diagram \ref{diagram}, keeping in
mind that  $\deg Z =  \chi \,\coo_Z$ and  $\deg_R \cae = \chi  \, \cae
-\chi \, \coo_R$ for an invertible sheaf $\cae$ on $R$.

\end{proof}

\bp\label{moduli}
Given a triple $\{Z,R,P \}$ of  closed subschemes of $F$ as above, the
set of curves  $C \subset X$ with $T(C)=\{Z,R,P \}  $ is in one-to-one
correspondence with an open subset of the vector space
$$
\mbox{\rm H}^{0} (R, \coo_{R} (Z+P-F)) \cong
\mbox{\rm Hom}_{R} (\coo_{R}(-P),  \coo_{R} (Z-F)  ).
$$
\ep

\begin{proof}
We study the fibres of the map $C \mapsto T(C)$.
We have seen that the bottom row of diagram (\ref{diagram})
is a sequence of $\coo_{R}$-modules, hence tensoring with $\coo_R$
we obtain a new diagram
\begin{equation}\label{phi}
\begin{CD}
0 @>>> \coo_{R}(-F)  @>{\tau}>>
\ideal{P} \otimes \coo_{R}
@>{\pi}>> \coo_{R} (-P) @>>> 0  \\
&& @V{=}VV  @V{\phi}VV @VVV \\
0 @>>> \coo_{R}(-F)  @>{\sigma}>> \coo_{R}(Z-F)
@>{\gamma}>> \coo_{Z}(-P) @>>> 0
\end{CD}
\end{equation}
in which $\phi$ is surjective, the top row is the conormal sequence of
$P$ in $X$ restricted to $R$, and the bottom row is obtained by dualizing
$   \exact{\ideal{Z,R}}{\coo_R}{\coo_Z}$.  It is clear that any
surjection ${\phi}$ with  ${\phi} \circ \tau = \sigma$  yields a curve
with triple $\{Z,R,P\}$. As in \cite{2h}, we obtain a
one-to-one correspondence between curves  $C$ in $X$ with triple $\{Z,R,P \}$
and surjections ${\phi}$ satisyfing ${\phi} \circ \tau = \sigma$.
Since the  cokernel of $\tau$ is isomorphic to $\coo_R (-P)$, the set of
such surjections may be identified with an open subset of
$\mbox{\rm  Hom}_{R} (\coo_{R}(-P), \coo_{R} (Z-F) )$.
\end{proof}

It is useful to know when a triple actually arises from a curve.

\bp \label{exist}
Let  $\{Z,R,P \}$  be   a  triple  of   subschemes  of  $F$   as  in
proposition ~\ref{1.1}. Suppose that
\begin{enumerate}
       \item\label{a} $\mbox{\rm H}^{1}   (R, \coo_{R}(Z+P-F)) =  0$; and
       \item\label{b} the map
       $\mbox{\rm    H}^{0} (\coo_{R} (Z+P-F)) \otimes \coo_{R}
      \rrr \coo_{Z}$ induced by $\gamma$ is
       surjective.
\end{enumerate}
Then the set of curves $C \subset X$ with
$T(C) = \{Z,R,P \}$ is parametrized by a non-empty open subset
$U \subset \mbox{\rm H}^{0}  (R, \coo_{R}  (Z+P-F))$ of dimension
$\deg Z  + \chi \coo_{R} (P-F)$.
\ep

\begin{proof}
The triple $\{Z,R,P \}$ gives rise to the two exact rows of
diagram~(\ref{phi}).
Condition $(1)$ gives
$$
\mbox{\rm Ext}^{1}(\coo_{R}(-P),\coo_{R}(Z-F)) \cong
\mbox{H}^{1}(\coo_{R}(Z+P-F))=0,
$$
hence there exists
$\phi_{0} \in \mbox{\rm Hom} (\ideal{P}  \otimes \coo_R, \coo_R (Z-F))$
such that $\phi_{0} \circ \tau = \sigma$. Moreover, any such
morphism $\phi$ can be written $\phi = \phi_{0} + \alpha \circ \pi$ for
$\alpha \in \mbox{\rm Hom} ( \coo_{R} (-P), \coo_R (Z-F)) \subset
\mbox{\rm Hom} ( \ideal{P} \otimes \coo_{R}, \coo_{R} (Z-F))$.

Let ${\overline {\phi_{0}}}:\coo_{R}(-P) \ra \coo_{Z}(-P)$ be the
morphism induced by $\phi_{0}$. The snake lemma shows that the morphism
$\phi_{0} + \alpha \circ \pi$ is surjective if and only if
${\overline {\phi_{0}}} + \gamma \circ \alpha$ is. Tensoring with
$\coo_{R}(P)$, we view $\alpha$ as a global section of
$\coo_{R}(Z+P-F)$. The images of these global section under
$\gamma$ generate $\coo_{Z}$ of by condition $(2)$.
Since $Z$ is finitely supported, it follows
that for a general such section $s \in H^{0} \coo_R (Z+P-F)$, the global
section $\gamma(s) + {\overline {\phi_{0}}}(1)$ is a unit in $\coo_{Z,z}$
at each point $z \in Z$. Thus $\alpha$ with $\alpha(1)=s$ corresponds
to a surjective morphism $\phi$.
\end{proof}

\bex
The hypotheses of Proposition~\ref{exist} are not necessary
for the existence of a curve $C$ with a given triple. Let
$F \subset \pso$ be a smooth surface of degree $d = \deg F$ which contains a
line $L$. The effective divisor $X = 2F$ on $\pso$ is a ribbon
over $F$ which contains all double lines $C$ supported on $L$. If
$p_{a}(C) \neq 1-d$, then $C \not \subset F$ and the triple $T(C)$
must take the form $\{Z,L,L\}$, where $Z$ is an effective
divisor of degree $d-1-p_{a}(C)$ on $L$. Since $\coo_{L}(Z+L-F) \cong
\omega_{L}(\deg Z + 4 - 2 \deg F)$ it is clear that
$H^{1}( \coo_{L} (Z+L-F)) \neq 0$ and
$H^{0}( \coo_{L} (Z+L-F))=0$ for $d > > 0$, hence neither hypothesis of
\ref{exist} hold, yet the existence of $C$ shows that the there are curves with
the triple $\{Z,L,L\}$. Replacing $L$ by any smooth curve gives
similar examples.
\eex

\br\label{practical} The following practical conditions
imply the hypotheses of Prop. \ref{exist}:
\begin{enumerate}
       \item $\mbox{H}^{1} (R, \coo_{R}(Z+P-F)) = 0$ and
       $\coo_{R}(Z+P-F)$ is generated by global sections.
       \item
       $\mbox{H}^1 (R, \coo_{R}  (Z+P-F-H))=0 $ for some very ample 
divisor $H$ on $R$.
       \item
       $\mbox{H}^{1} (R, \coo_{R}(P-F))=0$.
\end{enumerate}

Indeed, the first condition is clearly stronger than the hypotheses 
of~\ref{exist}.
The second condition implies the first by Mumford's regularity theorem.
If the third condition is satisfied,  then for {\em any} effective 
generalized divisor
       $Z \subset R$ the exact sequence
       $$0 \ra \coo_{R} (P-F) \ra \coo_{R} (Z+P-F)
       \stackrel{\gamma}{\ra} \coo_{Z} \ra 0$$
       shows that $\mbox{H}^{1} (R, \coo_{R} (Z+P-F))=0$ and that
       $\gamma$ is surjective on global sections, which implies hypothesis
       $(2)$ of~\ref{exist} since $Z$ has finite length.

\er

\br\label{split}
Perhaps the simplest situation occurs when the restriction to $R$ of the
conormal sequence associated to $P \subset F$ (the top row of diagram
\ref{phi}) splits. This happens in case (3) of Remark \ref{practical}.
\begin{enumerate}
       \item This splitting occurs if and only if
       the triple $\{\emptyset,R,P\}$ arises from a curve $C$, since in
       this case $\coo_{Z}=0$ and $\sigma$ is the identity map. In
       case $R=P$ is a general smooth space curve, we expect
       that $\{\emptyset,R,P\}$ does {\em not} arise from a curve $C$, as this
       would be equivalent to splitting of the normal bundle.

       \item If $P$ is the intersection of a surface $F \subset {\pp}^{n}$
       with a hypersurface $H$ of degree $d$, then restricting the
       natural map $\coo_{{\pp}^{n}}(-d) \ra \ideal{P}$ to $R$ provides
       such a splitting and the triple $\{\emptyset,P,P\}$ arises from
       the curve $X \cap H$. If $F$ is a general surface of
       degree $\geq 4$ in $\pso$, then every curve $P \subset F$ arises
       in this way since $\Pic F \cong \Z$ with
       $\coo_{S}(1)$ as generator \cite{gh}.

       \item When the splitting {\it does} occur, there is no
       obstruction to finding maps $\phi$ such that $\phi \circ \tau = \sigma$
       and the existence of a surjective such $\phi$ is {\em equivalent}
       to condition (2) of Proposition \ref{exist}.
\end{enumerate}
\er

\bex
If $X=2H$ is a double plane in $\pso$, then {\it every} triple
$\{Z,R,P\}$ with $Z$ Gorenstein arises from a curve $C \subset X$ by
Remark \ref{practical}(3) cf.~\cite{2h}.
\eex

In the next three examples, we consider behavior of triples
$\{Z,R,P\}$ for the double quadric $X=2Q \subset \pso$, using the
standard isomorphism $\Pic \; Q \cong \Z \oplus \Z$ \cite[II, 6.6.1]{AG}.

\bex\label{quadric1}
If $P$ has type $(a,b)$ with $a,b > 0$ and $R < P$, then {\it every}
triple $\{Z,R,P\}$ with $Z$
Gorenstein arises from a curve $C \subset X$. Indeed, the exact sequence
$$0 \ra \coo_{Q}(P-R-Q) \ra \coo_{Q}(P-Q) \ra \coo_{R}(P-Q) \ra 0$$
yields the vanishing of Remark \ref{practical}(3). Since Remark
\ref{split} applies here, it is common for the normal bundle of a curve
$P \subset Q$ to split when restricted to a {\em proper} subcurve $R$,
while this is quite rare when $R=P$ \cite{hulek}.
\eex

\bex\label{quadric2}
If $0 < a \leq b$ and $R=P$, then $\mbox{H}^{1}(\coo_{R}(P-Q)) \cong k$ and
condition \ref{practical}(3) fails.
\begin{enumerate}
     \item If $a=b$, then $P$ is a complete intersection and the
           triple $\{\emptyset,P,P\}$ arises from the complete intersection
      of $X$ and a surface of degree $a$ containing $P$.
      Not every triple $\{Z,P,P\}$ arises from a
      curve, however: if $P$ has type $(1,1)$ and $\deg (Z) =1$,
           the triple $\{Z,P,P\}$ could only be associated to a curve
           of degree $4$ and genus $2$ by (\ref{1.1}), but there is no such
           curve in $\pso$ \cite[3.1 and 3.3]{genus}.
      On the other hand, if $\deg (Z) \geq 2$ and $P$ is a smooth conic,
      then there exists a curve $C \subset X$ with $T(C) = \{Z,P,P\}$ by
      Proposition~\ref{exist}.

     \item If $1=a < b$ or $(a,b)=(1,0)$ and $P$ is a {\it smooth} rational
           curve, then the triple $\{Z,P,P\}$ arises from a curve if and only
      if $ \deg (Z) > 0$ (condition (2) of ~\ref{exist} fails when
      $\deg (Z) =1$, but any nonzero map $\phi$
      in diagram (\ref{phi}) is surjective in this case). The normal bundle
      splits (as described in \cite[Theorem 1]{hulek} over $\Cpl$), but the
      top row of diagram (\ref{phi}) does not.

     \item $1 < a < b$ and the pair $Z \subset P$ is sufficiently
           general with $\deg (Z) > 0$, then the triple $\{Z,P,P\}$ arises
      from a curve on $X$. Since the proof uses a degeneration argument,
      we postpone it until the following section on families (see Example
      \ref{degen} below).
\end{enumerate}
\eex

\bex\label{quadric3}
Now suppose that $P$ has type $(0,b)$ for some $b > 0$.
\begin{enumerate}
      \item Suppose that $R \subset P$ is a disjoint union of {\it reduced}
          lines. Then applying Example \ref{quadric2}(2)
     above to each line $L \subset R$, we see that the triple
     $\{Z,R,P\}$ arises from a curve $C \subset X$ if and only if
     $Z \cap L \neq \emptyset$ for each line $L \subset R$
     if and only if
     $\mbox{H}^{1}(\coo_{R}(Z+P-Q)) \cong \mbox{H}^{1}(\coo_{R}(Z-Q))=0$.

      \item Let $R \subset P$ be a double line on $Q$. In this case
          $Z$ need not be contained in the underlying reduced line.
     In fact, if $L$ is the underlying support, then the triple
     $Z \subset R \subset P$ satisfies the conditions of
     \ref{practical}(2) if $2 \leq \deg (Z \cap L) \leq \deg (Z) -2$.
          To see this, let $W = Z \cap L $ and let $Y$ be the residual scheme
          to $W$ in $Z$. Since $R^{2}=0$, the sequence relating $Y$ and $W$ to
          $Z$ takes the form
          $$ 0 \rrr \ideal{Y,L} \rrr \ideal{Z,R} \rrr \ideal{W,L} \rrr 0$$
          and applying $\Hom_{\coo_{R}}(-,\coo_{R})$ yields the exact sequence
          $$0 \rrr \coo_{L}(W) \rrr \coo_{R}(Z) \rrr \coo_{L}(Y) \rrr 0$$
          (one checks locally that $\Ext^{1}_{\coo_{R}}(\coo_{L},\coo_{R})=0$
          and $\Hom_{\coo_{R}}(\coo_{L}(a),\coo_{R}) \cong \coo_{L}(-a)$ by
          \cite[III,6.7]{AG}). Tensoring by $\coo_{R}(-Q-H)$ and taking the
          long exact cohomology sequence now gives the desired vanishing.
          One can formulate more complicated criteria for higher order
     multiple lines on $Q$.
\end{enumerate}
\eex

\section{Families}
In this section we study families of curves in $X$ and their
corresponding triples. We prove that, if $V$ denotes
the open subset of the flag Hilbert scheme $D$ consisting of
triples satisfying the conditions of Proposition \ref{exist},
the set of curves $C$ with $T(C) \in V$ is an open dense subset
of an affine fibre bundle over $V$ (\ref{structureofpi}).
Combining this with the structure of
the projection maps on Hilbert flag schemes (\ref{dominant}), we find that a
curve $C$ with triple $\{Z,R,P\}$ satisfying the first condition of
Proposition \ref{exist} is the flat limit of curves with triples
$\{Z^{\prime},R^{\prime},P^{\prime}\}$ for which $R^{\prime}$ is
general (\ref{limits}).

We first extend our constructions to the relative case.  Let
$Sch_{k}$  be the  category  of locally  Noetherian  schemes over  the
ground field $k$.
For $S \in Sch_{k}$, let $H(S)$ be the set of families of curves
$C \subset X \times S$ such that the sheaves $\coo_{C},
\coo_{C \cap (F  \times S)}$ and
$\cae_C = \Ext^{1}_{\coo_{F \times S}}
(\ideal{C \cap (F  \times S)},\coo_{F \times S})$
are all flat over $S$.
Then $H: Sch_{k} \rrr Sets$ defines a contravariant functor:
if $\phi:T \rrr S$ is a morphism in $Sch_{k}$, we define
$H(\phi): H(S) \ra H(T)$ by sending a family $C \in H(S)$ to its
pull-back $C_{T} = C \times_{S} T$. We have to check that this is well
defined, i.e., that $C_T \in H(T)$. Here the point is that
$\cae_{C_T}$ is the pullback of $\cae_{C}$: indeed,  on  the
fibres  we  have $\Ext^{2} (\ideal{C_{s}.F},\coo_{F})=0$, so the theorem of
base change for the $\Ext$ functors ~\cite{BPS,JS} tells us that $\cae_C$
commutes with base change - that is, the natural map
$({\rm Id}_{F} \times \phi)^{*} \cae_C \ra \Ext^{1}_{\coo_{F \times T}}
(\ideal{C_{T} \cap (T \times F)},\coo_{F \times T})$
is an isomorphism.

We now claim that to a family  of curves $C \in H(S)$ we can associate
a  triple $T(C) =  \{Z,R,P\}$ where  $Z \subseteq  R \subseteq  P$ are
closed subschemes of $F \times S$,  flat over $S$, such that for every
closed  point  $s  \in  S$  the  triple  $\{Z_{s},Y_{s},P_{s}  \}$  is
precisely the triple $T(C_{s})$.

To construct $P$, we need to show that the sheaf
$$
\mathcal{H}_C =
\Hom_{\coo_{F \times  S}}  (\ideal{C.F},\coo_{F \times
S})
$$
is an  invertible sheaf on $F  \times S$. By definition  of $H(S)$, we
know that $\cae_C $ is flat
over $S$, and its formation commutes with base change.  The theorem of
base change  for the  functors  $\Ext^{i}$ implies that  $ \mathcal{H}_C $
itself is flat over $S$ and commutes with base change.

In particular, the natural map
$$
\mathcal{H} \otimes k(s) \rrr \Hom (\ideal{C_{s}.F},\coo_{F})
$$
is  an isomorphism  for  every closed  point  $s \in  S$.  Thus  the
restriction  of $\mathcal{H}$ to  each fibre is an  invertible sheaf,
hence so is $\mathcal{H}$.

By a  standard argument \cite[7.4.1]{kollar},  the inclusion
$\ideal{C.F}  \hookrightarrow  \coo_{F  \times  S}$ defines  a  global
section of $\mathcal{H}$ whose zero scheme is an effective Cartier
divisor $P \subset F \times S$, flat over $S$.

Now define $Z (C) \subset F \times S$ to be the residual scheme to $P$ in
$C \cap (F \times S)$, so that
$\ideal{C \cap (F \times S)}=  \ideal{P}\ideal{Z}$.  To see that $Z$ is
flat  over $S$  we note  $\coo_{Z} (-P)  \cong \ideal{P,C \cap (F \times S)}$
and use \cite[7.4.1]{kollar}.

Finally, define $R(C) \subset X \times S$ to be the residual
scheme to the intersection of $C$ with $F \times S$. The exact sequence
$$
   0 \rrr
\coo_{R}(-F \times S) \rrr \coo_{C} \rrr
\coo_{C.F}  \rrr 0
$$
shows that $R$ is flat over $S$, and that for each $s \in S$ the fibre
$R_{s}$ is the residual scheme to the intersection of $C_{s}$ with
$F$.  Since $Z_{s} \subseteq R_{s} \subseteq P_{s}$ for
each $s \in S$, we have $Z \subseteq R \subseteq P$.

Summing up,  to any  $C \in H(S)$  we can  associate a triple  $T(C) =
\{Z,R,P \}$ where  $Z \subseteq Y \subseteq P$,  are closed subschemes
of $F \times  S$, flat over $S$, and this construction is compatible with base
change. Thus we have  a natural  transformation $T: H  \ra D$
where  $D$ is the functor that to a  scheme $S$ associates flags
$Z \subset R \subset P \subset F \times S$, with $Z$, $R$, $P$ flat over
$S$, $Z$ zero dimensional, and $R \subset P$ effective Cartier divisors.

Both $H$ and  $D$ are represented by quasiprojective  schemes. This is
well known for $D$.  Using Mumford's flattening stratification, we see
$H$ is representable by a subscheme of the Hilbert scheme of curves in
$X$. Since giving the  Hilbert  polynomials of  $C$,
$C \cap (F \times S)$  and
$\cae_C$ is the same as  giving the Hilbert polynomials of $Z$, $R$ and
$P$, $H$  is represented by the disjoint  union of locally closed
subschemes $H_{z,r,p}$ of  the Hilbert scheme of curves  in $X$, where
$\{z,r,p\}$  vary  in the  set  of  possible  Hilbert polynomials  for
$Z$, $R$ and $P$.  Furthermore,  the natural transformation $T$ induces
a morphism of schemes $t: H_{z,r,p} \ra D_{z,r,p}$.

\bt \label{structureofpi}
Let $V \subset D_{z,r,p}$ be the open subscheme corresponding to
triples $\{Z,R,P\}$ satisfying $H^{1}(\coo_{R}(Z+P-F))=0$. Then
the map $t^{-1}(V) \ra V$ has the structure of an open immersion
followed by a projection from an affine bundle over $V$.
\et

\begin{proof}
Given a triple $\{Z,R,P\} \in D(S)$, we define
$$
\coo_R (Z - F \times S) =
{\mathcal H}om_{\coo_R} ( \ideal{Z,R} , \coo_R (-F \times S)).
$$
If  $s \in S$  is a  closed point,  we have  $\cae xt^{1}_{\coo_{R_s}}
(\ideal{Z_s,R_s},   \coo_{R_s}  (-F)   )   =  0$   because  $R_s$   is
Gorenstein. It follows that
$\coo_R (Z - F \times S)$
is flat over $S$ and its formation
commutes with base change \cite{BPS,JS}, so that for every
morphism $g:  T \ra S$ in  $Sch_{k}$ the pull  back of $\coo_R (Z  - F
\times S)$ is $\coo_{R_T} (Z_T - F \times T)$.
Hence there is a functor
$A=A_{z,r,p}$ that assigns to the  scheme $S$ the set of flat families
of  flags $Z  \subset R  \subset P  \subset F  \times S$  with Hilbert
polynomials $z$,$r$,$p$ along with a morphism
$\phi:\ideal{P} \otimes \coo_{R} \ra \coo_{R} (Z-F \times S)$.

The exact sequence
\begin{equation} \label{m}
0  \ra \coo_{R}  (- F  \times S)  \stackrel{\tau}{\ra} \ideal{P}
\otimes \coo_{R}
\stackrel{\pi}{\ra} \coo_{R}  (- P) \ra 0.
\end{equation}
and the sequence
\begin{equation} \label{l}
0 \ra
\coo_{R} (- F \times S) \stackrel{\sigma}{\ra} \coo_R (Z - F \times S)
\ra
\cae xt^{1}(\coo_{Z_S} , \coo_{R_S} (- F \times S) ) \ra 0
\end{equation}
obtained by dualizing
$$\exact{\ideal{Z,R}}{\coo_{R}}{\coo_{Z}}$$
are both compatible with base change, thus $A$ has a subfunctor
$M=M_{z,r,p}$ corresponding to
morphisms $\phi$ satisfying $\phi \circ \tau = \sigma$.

Now we claim  that $H_{z,r,p}$ is an open subfunctor of $M$.
Indeed, given $C \in H(S)$, we may write a diagram analogous to diagram
(\ref{diagram}):
\begin{equation} \label{fourth}
\begin{CD}
&& && && 0&& \\
&&&& && @VVV \\
&&&& \ideal{C} \otimes \coo_R @>>>
\ideal{Z,R} (-P)@>>> 0 \\
&&&&  @VVV @VVV \\
0 @>>> \coo_{R} (-F \times S) @>\tau>> \ideal{P} \otimes \coo_R @>>>
\coo_{R} (-P) @>>> 0\\
&&@VVV   @VV{\phi}V  @VVV \\
0 @>>> \coo_{R} (-F \times S) @>\sigma>> \scrl @>>>  \coo_{Z} (-P) @>>> 0    \\
&&&& @VVV @VVV \\
&&&& 0 && 0&& \\
\end{CD}
\end{equation}

As in the proof of~\ref{1.1}, we obtain a morphism
$\psi: \scrl \ra \coo_R (Z-F \times S)$. These sheaves are flat over
$S$ and compatible with pull back. Since $\psi$ induces isomorphisms
$\psi_s$ on the fibres by the proof of~\ref{1.1},
$\psi$ is an isomorhism. Thus the diagram gives us a morphism
$\phi:\ideal{P} \otimes \coo_{R} \ra \coo_{R} (Z-F)$ with
$\phi \circ \tau = \sigma$, and we obtain a natural transformation from
$H$ to $M$ that makes $H$ into a subfunctor of $M$. It is open because
it corresponds to the open condition that the map $\phi$ be surjective.


It remains to show that when we take inverse images over $V \subset D$,
the induced map $M_{V} \stackrel{t}{\ra} V$ has the structure of an
affine bundle.
       Let $U \subset V$ be an affine open set equipped with
       universal flat flag
       $$\begin{matrix}
       Z  & \subset & R & \subset & P & \subset & F \times U\\
       &&&&&& \downarrow f \\
       &&&&&& U.\\
       \end{matrix}$$
       Since $\mbox{H}^{1} (\coo_{R_{u}}(Z_{u}+P_{u}-F))=0$ for each $u \in U$,
       we deduce \cite[III, 8.5 and 12.9]{AG} that $R^{1} f_{*}
       \coo_{R}(Z+P-F)=0$ and hence that
       $${\rm {Ext}}^{1}(\coo_{R}(-P),\coo_{R}(Z-F)) \cong
       H^{1}(\coo_{R}(Z+P-F))=0.$$
       In particular, there exists $\phi_{0}:\ideal{P} \otimes \coo_{R}
       \ra \coo_{R}(Z-F)$ such that $\phi_{0} \circ \tau = \sigma$.

Now let $G: Sch_{U} \ra Sets$ be the functor that to a scheme $T$
over $U$ associates
the set
$$ G(T)= \mbox{Hom}_{R_T} ( \coo_{R_T} (-P_T),\coo_{R_T} (Z_T - F
\times T).$$
By the lemma \ref{bundle} below,
       $\cae = f_{*} \cah om_{\coo_{R}} (\coo_{R}(-P), \coo_{R}(Z-F))$ is
       locally free  on $U$, and $G$ is represented by the geometric vector
       bundle $B \stackrel{p}{\ra} U$ whose sheaf of sections is $\cae$.
       Thus there is a universal map
       $\alpha:\coo_{R_{B}}(-P_{B}) \ra \coo_{R_{B}}(Z_{B}-F)$
       on the pullback of the universal flag to $B$. We now show that
       the pair $(B,\phi=p^{*}(\phi_{0})+\alpha \circ \pi)$ represents
       $M_{U}$.

       To this end, let $S$ be a scheme,
       $Z_{S} \subset R_{S} \subset P_{S} \subset F \times S$ be a flag
       corresponding to a map $h: S \ra D$ that factors through $U$, and
       $\psi:\ideal{P_{S}} \otimes \coo_{R_{S}} \ra \coo_{R_{S}} (Z-F)$
       be a map satisfying
       $\psi \circ \tau_{S} = \sigma_{S}$.  By construction the map
       $\psi - h^{*}(\phi_{0})$ is the image of a map
       in $\mbox{Hom} (\coo_{R_{S}} (-P_{S}), \coo_{R_{S}} (Z_{S}-F_{S}))$,
       hence the universal property of $B \ra S$ yields a unique lifting
       ${\tilde h}:S \ra B$ of $h$. Moreover, it is clear from
       construction that $\psi = {\tilde h}^{*}(\phi)$. This shows that
       $(B,\phi)$ represents $M_{U}$, finishing the proof.

\end{proof}
The following lemma, which we used in the above proof, is an immediate
consequence of the theorems of base change
for cohomology and for the $\Ext$ functors.
\bl\label{bundle}
Let $f: R  \rightarrow U$ be a morphism  of locally Noetherian schemes
over  $k$,   and  let $\caf$, $\cag$ be  coherent  sheaf   on  $R$.   Let
$G=G_{\caf,\cag}:  Sch_{U}  \rightarrow  Sets$  be  the  contravariant
functor that  to a locally  Noetherian $U$-scheme $T$ associates the set
$$
G(T) = \mbox{\em Hom}_{R_T} ( \caf_T,\cag_T)
$$
where $R_T$, $\caf_T$, $\cag_T$ are the base extensions to $T$.  Suppose
that $f$ is projective and flat, and $\caf$,$\cag$ are flat over $U$.
Furthermore, suppose that for every  point $u \in U$:
\begin{enumerate}
\item
$\cae xt^{1}_{\coo_{R_u}} ( \caf_u, \cag_u ) = 0$;
\item
$\text{\em H}^{1} (R_u, \cah om_{\coo_{R_u}} ( \caf_u, \cag_u ) ) = 0$.
\end{enumerate}
Then the sheaf $\, \cae =  f_{*} \cah om_{\coo_{R}} ( \caf, \cag )$ is
locally free  on $U$, and $G$  is represented by  the geometric vector
bundle over $U$ whose sheaf of sections is $\cae$.
\el

\bc \label{comps1}
Let $Y$ be an irreducible component of  $D_{z,r,p}$ and let $U \subset Y$ be
the open subset consisting of triples $\{Z,R,P\}$ for which
$H^{1}(\coo_{R}(Z+P-F))=0$. If $t^{-1}(U)$ is nonempty,
then ${\overline {t^{-1}(U)}}$ is an irreducible component of $H_{z,r,p}$.
\ec

\begin{proof}
       From the structure of $t$ given in Theorem \ref{structureofpi},
       ${t^{-1}(U)} \subset H_{z,r,p}$ is an irreducible open subset of
       ${t^{-1}(Y)}$.
       Let $W$ be an irreducible component of $H_{z,r,p}$ containing
       ${\overline {t^{-1}(U)}}$. Then $t(W)$ is irreducible and contains
       a nonempty open subset of $U$ ($t|_{t^{-1}(U)}$ is an open map
       by \ref{structureofpi}),
       hence $Y= {\overline U} \subset \overline{t(W)}$.
       Since $Y$ is an irreducible component, we must have
       $Y=\overline{t(W)}$, hence $W \subset {t^{-1}(Y)}$.
       It follows that $t^{-1}(U) = t^{-1}(U) \cap W$ is a nonempty open
       subset of $W$ and $W = {\overline {t^{-1}(U)}}$.
\end{proof}

\br
Note that $t^{-1}(V)$ (resp. $t^{-1}(U)$) may be empty in Theorem
\ref{structureofpi} (resp. Cor. \ref{comps1}), as in the case of the smooth
conic of type $(1,1)$ on the quadric surface and $\deg Z = 1$ (Example
\ref{quadric2}(1)). These sets are guaranteed to be nonempty if
there exists a triple $\{Z,R,P\}$ in $V$ (resp. $U$) satisfying
condition two of Proposition \ref{exist}.
\er

To use Corollary~\ref{comps1}, we need to understand the Hilbert scheme of
flags $D_{z,r,p}(F)$. Now $D_{z,r,p}$
breaks up as the disjoint union of closed subschemes $D_{z,\xi,\eta}$  where
$\xi$ (resp. $\eta$) varies in the set of numerical equivalence
classes of divisors
in $F$ with Hilbert polynomial $r$ (resp. $p$).
We have a decomposition
$$D_{z,\xi,\eta} \cong D_{z,\xi} \times H_{\eta - \xi}$$
where $D_{z,\xi}$ denotes the Hilbert scheme of flags $Z \subset R \subset F$,
with $Z$ zero dimensional of degree $z$, and $R$ an effective divisor of class
$\xi$, and $H_{\eta - \xi}$ is the Hilbert scheme of effective divisors in $F$
of class $\eta - \xi$ - this because we can tack on the effective divisor $P-R$
after choosing the flag $Z \subset R$. The following lemma helps to
identify the irreducible components of the Hilbert flag scheme.

\bl\label{dominant}
Let $q:D_{z,\xi} \ra H_{\xi}$ be the projection.
Then $q$ is surjective and maps generic points of $D_{z,\xi}$ to
generic points of $H_{\xi}$.
\el

\begin{proof}
The argument is due to Brun and Hirschowitz \cite[3.2]{brun-hirsch}.
Since $q$ is proper and surjective, it is enough to show that, if
$A$  is an irreducible component of $ D_{z,\xi}$ and $B$ is an
irreducible component of
$H_{\xi}$ that contains  $q(A)$, then $B= q(A)$.
Let $M = \Hilb_{z}(F)$. $D_{z,\xi}$ is constructed as the scheme of zeros
of a global section of a rank $z$ vector bundle on $M \times H_{\xi}$
\cite{kleppe,sernesi}.

Thus the codimension of $D_{z,\xi}$ in $M \times H_{\xi}$ is $\leq z$
at each point.
       In particular, the irreducible component $A$ has
       dimension at least $\dim B + z$.

       On the other hand, let $J \subset B$ denote the image
       of $A$. The fibre over any fixed curve $Y \in B$ has dimension
       $\leq z$ by the theorem of Brian\c con \cite{briancon,iar} which
       describes the punctual Hilbert scheme. It follows that
       $$\dim B + z \leq \dim A \leq \dim J + z,$$
       hence these are equalities and $J=B$.

\end{proof}

\br
If Z is Cartier on R, it follows from deformation theory
that the map q of lemma 3.4 is smooth at the point $(Z,R)$ of $D$,
because $H^1(\cnn_{Z,R})=0$.
\er

\br\label{irred}
If $B \subset H_{\xi}$ is an irreducible component whose general
member is a smooth connected curve, then $q^{-1}(B)$ is an
irreducible component of $D_{z,\xi}$ (\cite[4.3]{2h}). Indeed, the
irreducible components of $D_{z,\xi}$ contained in $q^{-1}(B)$ map
dominantly to $B$, but the general fibre of $q$ is irreducible, so there
is only one such component.
\er

\bex\label{ex1}
(1) If $F = \Ptwo$, the class of a divisor is determined by its degree $d$.
If $B={\PP} H^{0}(\coo_{\Ptwo}(d))$, then $D_{z,d} = q^{-1}(B)$ is
irreducible by Remark \ref{irred}. It now follows from
Corollary \ref{comps1} and Remark \ref{practical}(3) that the schemes
$H_{z,r,p}$ are irreducible. In fact, their closures are precisely the
irreducible components of $H_{d,g} (X)$ \cite[5.1]{2h}.
\eex

\bex\label{ex2}
If $F=Q \subset \pso$ is the smooth quadric surface,
then the numerical equivalence class of a divisor is determined
by its bidegree $(a,b)$ \cite[II,6.6.1]{AG} and $H_{(a,b)} = |\coo_{Q}(a,b)|$
is a projective space.
If $a$ and $b$ are both positive, then the general element of $H_{a,b}$
is smooth and irreducible, hence $D_{z,(a,b)}$ is irreducible by
Remark \ref{irred}.

If $a=0$ and $b > 0$, then the general element in $H_{(a,b)}$ is a
disjoint union of $b$ lines and $D_{z,(a,b)}$ has irreducible components
corresponding to various partitions of $z$ as a sum of $b$
non-negative integers, depending on how the zero-dimensional scheme
$Z$ is distributed among the generic lines in the family. In
particular, $q^{-1}(H_{(a,b)})$ is not irreducible unless $z \leq 1$
or $b=1$.
\eex

We now prove that if $T(C)=\{Z,R,P\}$ satisfies
$H^{1}(\coo_{R}(Z+P-F))=0$, then a general deformation
of $R$ lifts to a deformation of $C$.

\bt\label{limits}
Let $C \subset X$ be a curve with triple $T(C)=\{Z,R,P\}$
such that $H^{1}(\coo_{R}(Z+P-F))=0$.
Suppose that $B$ is an irreducible
component of $H_{r} (F)$ containing $R$. Then there is an irreducible
component $W$ of $H_{z,r,p}$ containing $C$ such that the natural
map $H_{z,r,p} \ra H_{r}$ induces a dominant map $W \ra B$.
\et

\begin{proof}
By Lemma~\ref{dominant} there is an irreducible component
$X \subset D_{z,r}$ containing $(Z,R)$ such that
$q(X)=B$. Since
$D_{z,\xi,\eta} \cong D_{z,\xi} \times H_{\eta-\xi}$ (here $\eta$ and $\xi$
are the numerical equivalence classes of $P$ and $R$; see discussion
following Cor. \ref{comps1}), we obtain an
irreducible component $Y = X \times K$ of $D_{z,\xi,\eta}$ containing
$T(C)$ and mapping dominantly to $B$ for a suitable irreducible
component $K \subset H_{\eta-\xi}$. Letting $U \subset Y$ be the
open set of triples for which the vanishing occurs, $U$ is dense in $Y$
and the generic point of $U$ maps to the generic point of $B$. By
Lemma~\ref{comps1} $W=\overline{t^{-1}(U)}$ is an irreducible
component of $H_{z,r,p}$, and by construction the generic
point of $t^{-1} (U)$ maps to the generic point of $B$.
\end{proof}

\bex\label{thick}
The conclusion of Theorem \ref{limits} fails for a general thick
$4$-line $C$ of genus $g$ on the double quadric $X=2Q$ in $\pso$.
Recall that a {\it thick} $4$-line is a curve of degree $4$ supported
on a line $L$ and containing the first infinitesimal neighborhood $L^{(2)}$
\cite{banica}.
We claim that such a curve is {\bf not} a flat limit of disjoint unions
of double lines on $X$.
To see this, we first note that the family of
double lines of genus $g_{1}$ with fixed support is irreducible of
dimension $1-2g_{1}$ by \cite[1.6]{nthree}. Since the lines on $Q$
form a one-dimensional family, the disjoint unions of two double lines
of genera $g_{1}$ and $g_{2}$ form a family of dimension
$4-2g_{1}-2g_{2}=2-2g$.

On the other hand, the thick $4$-lines on fixed support $L$ are determined
by surjections in
$${\rm {Hom}} (\ideal{L^{(2)}}, \coo_{L}(-g-1)) \cong
{\rm {Hom}} (\coo_{L}(-2)^{3}, \coo_{L}(-g-1)) \cong 
H^{0}(\coo_{L}(-g+1)^{3})$$
by \cite[$\S 4$]{banica}, hence these form an irreducible family of
dimension $5-3g$. We are interested in the subset of
those which send the equation of $X$ to zero. If $L = \{x=y=0\}$
and $Q = \{xz-yw=0\}$, then $X = \{x^{2}z^{2}-2xyzw+y^{2}w^{2}=0\}$
and hence the thick $4$-lines with support $L$ lying on $X$ correspond
to the triples
$\{(a,b,c) \in H^{0}(\coo_{L}(-g+1)^{3}) : az^{2}-2zwb+cw^{2}=0\}$.
These form a vector subspace of codimension $-g+4$
(provided char $k \neq 2$), hence the family has dimension $1-2g$.
Varying the support line $L$ on $Q$, we obtain a family of dimension $2-2g$
and conclude that the general thick $4$-line $C$ cannot be the limit
of a family whose general member is a disjoint union of two double lines.
\eex

\bex\label{degen} In Example \ref{quadric2}(3) we claimed that if
$X = 2Q \subset \pso$ is the double quadric, then the general triple
$\{Z,P,P\}$ arises from a curve on $X$ if $\deg Z > 0$ and $P$ has
type $(a,b)$ with $1 < a < b$. We now explain why.

Let $C$ be a smooth rational curve of type $(1,b-a+1)$ on $Q$ and
$Z \subset C$ a divisor with $\deg Z > 0$.
By Example \ref{quadric2}(2), $H^{1} (\coo_{C} (Z+C-Q))=0$ and the triple
$\{Z,P,P\}$ arises from a curve ${\tilde C}$ on $X$. If $H$ is a
general hypersurface of degree $a-1$, then $H \cap {\tilde C}$
consists of $(a-1)(b-a+2)$ double points and $H \cap Z = \emptyset$.
Letting ${\tilde E}=H \cap X$, we have $T({\tilde E})=\{\emptyset,E,E\}$
where $E$ is a divisor on $Q$ of type $(a-1,a-1)$ (see \ref{quadric2}(1)).
The triple for ${\tilde C} \cup {\tilde E}$ has form
$\{{\tilde Z},C \cup E,C \cup E\}$, but $Z \subset {\tilde Z}$ by
local considerations and the genus formula forces $Z = {\tilde Z}$,
hence $T({\tilde C} \cup {\tilde E})=\{Z,C \cup E,C \cup E\}$.

Since $H^{1} (\coo_{C} (Z+C-Q))=0$, when we tensor the exact sequence
$$0 \ra \coo_{C}(C) \ra \coo_{C \cup E}(C+E) \ra \coo_{E}(C+E) \ra 0$$
by $\coo_{C \cup E}(Z-Q)$ we see that
$H^{1} (\coo_{C \cup E} (Z+C+E-Q))=H^{1} (\coo_{P} (Z+P-Q))=0$. Note here that
$H^{1} (\coo_{E}(Z+C+E-Q))=0$ via the exact sequence
$$0 \ra \coo_{E}(Z+E-Q) \ra \coo_{E}(Z+E+C-Q) \ra \coo_{E \cap 
C}(Z+E+C-Q) \ra 0.$$
We can now apply Theorem \ref{limits} and its proof to
${\tilde C} \cup {\tilde E}$. By Remark \ref{irred}(b), $B = H_{(a,b)}$ is
irreducible as is $D_{z,r} \cong D_{z,r,p}$, hence the
general triple $\{Z,P,P\}$ with $\deg Z > 0$ and $P$ of type $(a,b)$
arises from a curve.
\eex

\bex\label{3.3}
Let $W$ be a quasi-primitive triple line of type $(0,b)$ in $\pso$ for
some $b \geq 0$.
Then the underlying double line $D$ necessarily lies on a smooth
quadric surface $Q$ \cite[1.5]{nthree} and hence $W$ lies on the
double quadric $X=2Q$.
The associated triple is $T(W)=\{Z,L,D\}$, where $L$ is the support of $W$
and $Z \subset L$ is a divisor of degree $b+2$ by the genus formula
of Proposition \ref{1.1} ($g(W)=-2-b$ by \cite[2.3a]{nthree}).
If $H$ denotes the hyperplane divisor, then
$$\mbox{H}^{1}(\coo_{L}(Z+D-Q-H)) \cong \mbox{H}^{1} (\coo_{L}(b-1))=0$$
since $b \geq 0$ and Remark \ref{practical}(2) applies.
We deduce from Theorem \ref{limits} that $W$ is the
limit of a family of curves on $2Q$ whose general member is the
disjoint union of a line and a double line. This generalizes the
deformation used in the proof of \cite[3.3]{nthree}.
\eex

\bex\label{4lines}
This is the example that inspired the present paper. Let $R=P$ be a
double line $2L$ on the smooth quadric surface $Q \subset \pso$.
Let $c \geq b \geq 0$ be integers and let $Z \subset R$ be a divisor
consisting of $c-b$ simple points and $b+2$ double points which are
not contained in $L$. One can show \cite[3.2]{NS} that the triple
$\{Z,R,P\}$ arises from a general quasiprimitive 4-line $C$ of type
$(0,b,c)$. Since $Z$ contains $ \geq 2$ double points not
contained in $L$, the sufficient condition of Example
\ref{quadric2}(b) holds - thus the conditions of Proposition
\ref{exist} hold for the triple $\{Z,R,P\}$ and
we can apply Theorem \ref{limits} to see that $C$ is in the closure
of some family of disjoint unions of double lines, since the general member of
$|\coo_{Q}(0,2)|$ is a disjoint union of two lines.
\eex

\end{document}